\documentclass[12pt]{article}
\usepackage{amssymb}
\usepackage{latexsym,bm}
\usepackage{graphicx}
\usepackage{amsmath}
\usepackage{mathrsfs}

\date{}
\begin{document}
\title{Regular homogeneously traceable nonhamiltonian graphs\footnote{E-mail addresses:
{\tt huyanan530@163.com}(Y.Hu),
{\tt zhan@math.ecnu.edu.cn}(X.Zhan).}}
\author{\hskip -10mm Yanan Hu and Xingzhi Zhan\thanks{Corresponding author.}\\
{\hskip -10mm \small Department of Mathematics, East China Normal University, Shanghai 200241, China}}\maketitle
\begin{abstract}
A graph is called homogeneously traceable if every vertex is an endpoint of a Hamilton path. In 1979 Chartrand, Gould and Kapoor proved  that for every integer
$n\ge 9,$ there exists a homogeneously traceable nonhamiltonian graph of order $n.$ The graphs they constructed are irregular. Thus it is natural to consider the existence problem of regular homogeneously traceable nonhamiltonian graphs. We prove two results: (1) For every even integer $n\ge 10,$ there exists a cubic homogeneously traceable nonhamiltonian graph of order $n;$ (2) for every integer $p\ge 18,$ there exists a $4$-regular homogeneously traceable graph of order $p$ and circumference $p-4.$ Unsolved problems are posed.
\end{abstract}

{\bf Key words.} Homogeneously traceable; regular graph; circumference

{\bf Mathematics Subject Classification.} 05C38, 05C45, 05C76

\section{Introduction}

We consider finite simple graphs. The {\it order} of a graph is its number of vertices, and the {\it size} is its number of edges.  We denote by $V(G)$  the vertex set of a graph $G.$  The following concept is introduced by Skupie\'{n} in 1975 (see [3, p.185], and [4]). Note that the preprint of the 1984 paper [4]
was cited by the 1979 paper [2]. 

{\bf Definition 1.} A graph $G$ is said to be {\it homogeneously traceable} if every vertex of $G$ is an endpoint of a Hamilton path.

Obviously, hamiltonian graphs and hypohamiltonian graphs are homogeneously traceable. 
Chartrand, Gould and Kapoor [2] proved that for every integer $n$ with $3\le n\le 8,$ any homogeneously traceable graph of order $n$ is hamiltonian
and that for $n\ge 9,$ there exists a homogeneously traceable nonhamiltonian graph of order $n.$ This result was rediscovered in [1] where the term
``homogeneously traceable" was called ``fully strung". The homogeneously traceable nonhamiltonian graphs constructed in [2] are irregular
while the homogeneously traceable nonhamiltonian graphs constructed in [1] are also irregular except the Petersen graph of order $10$ which is cubic 
(i.e., $3$-regular). 
Thus it is natural to consider the existence problem of regular homogeneously traceable nonhamiltonian graphs.

In Section 2 we construct regular homogeneously traceable nonhamiltonian graphs, and in Section 3 we pose two unsolved problems.

\section{Regular homogeneously traceable nonhamiltonian graphs}

Given a vertex $v$  in a graph,  a {\it $v$-path} is a path with $v$ as an endpoint. We use $K_d$ to denote the complete graph of order $d,$ and use $N(v)$ to denote the neighborhood of a vertex $v.$ The notation ${\rm circum}(G)$ means the circumference of a graph $G.$

{\bf Definition 2.}  Let $v$ be a vertex of degree $d$ in a graph. {\it Blowing up $v$ into the complete graph $K_d$} is the operation of replacing $v$ by $K_d$ and adding $d$ edges joining the vertices of $K_d$ to the vertices in $N(v)$  such that the new edges form a matching.

The operation of blowing up a vertex of degree $4$ into $K_4$ is depicted in Figure 1.
\vskip 3mm
\par
 \centerline{\includegraphics[width=3.8in]{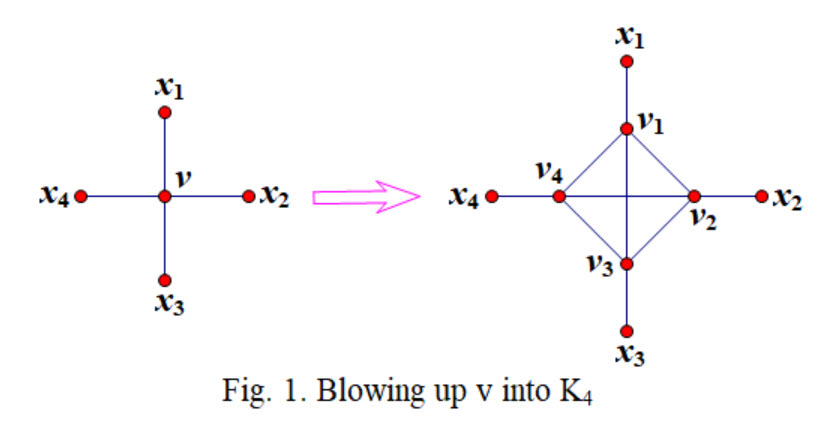}}
\par

{\bf Definition 3.} A graph $G$ is called {\it doubly homogeneously traceable} if for any vertex $v$ of $G,$ there are two Hamilton $v$-paths $P$ and $Q$
such that the two edges incident to $v$ on $P$ and $Q$ are distinct.

We will need the following two lemmas.

{\bf Lemma 1.} {\it Let $v$ be a vertex of degree $3$ in a doubly homogeneously traceable graph $G$ of order $n$ and circumference $c.$ Suppose $G^{\prime}$ is the graph obtained from $G$ by blowing up $v$ into $K_3.$ Then $G^{\prime}$ is also doubly homogeneously traceable. If $v$ lies in a longest cycle of $G,$  then $G^{\prime}$ has circumference $c+2.$  }

{\bf Proof.} Let $N(v)=\{x_1,x_2,x_3\}$ and suppose $v$ is blown up into $K_3$ whose vertices are $v_1,v_2,v_3$ such that $v_i$ is adjacent to $x_i$ for $i=1,2,3.$ Let $u\in V(G^{\prime}).$ If $u\not\in \{v_1,v_2,v_3\},$ there exist two Hamilton $u$-paths $P: u,\ldots,x_i,v,x_j,\ldots$ and $Q: u,\ldots,x_s,v,x_t,\ldots$ of $G$ where the two edges incident to $u$ on $P$ and $Q$ are distinct. Then $G^{\prime}$ has two Hamilton $u$-paths $u,\ldots,x_i,v_i,v_f,v_j,x_j,\ldots$ and $u,\ldots,x_s,v_s,v_g,v_t,x_t,\ldots$ where the two edges incident to $u$ are distinct.

Next suppose $u\in \{v_1,v_2,v_3\}.$ Without loss of generality suppose $G$ has two Hamilton $v$-paths $v,x_1,\ldots$ and $v,x_2,\ldots.$ Then $G^{\prime}$ has two Hamilton $v_3$-paths: $v_3,v_2,v_1,x_1,\ldots$ and $v_3,v_1,v_2,x_2,\ldots.$ Since $G$ is doubly homogeneously traceable, $G$ has a Hamilton $x_1$-path
$x_1,y,\ldots,x_i,v,x_j,$ $\ldots$ with $y\not=v.$ It follows that $G^{\prime}$ has two Hamilton $v_1$-paths: $v_1,x_1,y,\ldots,x_i,v_i,v_j,x_j,\ldots$
and $v_1,v_3,v_2,x_2,\ldots,$ where the two edges $v_1x_1$ and $v_1v_3$ are distinct. Similarly we can show that $G^{\prime}$ has two Hamilton $v_2$-paths
where the two edges incident to $v_2$ are distinct. This completes the proof that $G^{\prime}$ is doubly homogeneously traceable.

Now suppose $v$ lies in a longest cycle of $G.$  Let $\ldots,x_i,v,x_j,\ldots$ be a cycle of $G$ with length $c.$
Then $G^{\prime}$ contains the cycle $\ldots,x_i,v_i,v_f,v_j,x_j,\ldots$ which has length $c+2.$ Thus ${\rm circum}(G^{\prime})\ge c+2.$ On the other hand,
let $C$ be a longest cycle of $G^{\prime},$  which has length at least $c+2.$  If $C$ contains no vertex from the set $S=\{v_1,v_2,v_3\},$ it is also a cycle in $G$ and hence has length at most $c,$ a contradiction. Observe that every vertex in $S$ has exactly one neighbor outside $S.$
If $C$ contains a vertex in $S,$ then $C$ contains at least two vertices in $S.$ Note that the vertices in $V(C)\cap S$ appear consecutively on $C.$ Since $S$ is a clique and $C$ is a longest cycle, we deduce that $|V(C)\cap S|=3.$
Thus $v_rv_sv_t$ is a path  on $C$ with $\{r,s,t\}=\{1,2,3\}.$ Replacing this path by the vertex $v$ we obtain a cycle in $G$ which has length at most $c.$
Hence $C$ has length at most $c+2,$ implying that ${\rm circum}(G^{\prime})\le c+2.$ Combining this inequality with ${\rm circum}(G^{\prime})\ge c+2$ we obtain ${\rm circum}(G^{\prime})=c+2.$ \hfill $\Box$

{\bf Lemma 2.} {\it Let $v$ be a vertex of degree $4$ in a doubly homogeneously traceable graph $G$ of order $n$ and circumference $c.$ Suppose $G^{\prime}$ is the graph obtained from $G$ by blowing up $v$ into $K_4.$ Then $G^{\prime}$ is also doubly homogeneously traceable. If $v$ lies in a longest cycle of $G$ and $v$ lies in a clique
of cardinality $4,$ then $G^{\prime}$ has circumference $c+3$ and $G^{\prime}$ contains a vertex $v^{\prime}$ that lies in a longest cycle of $G^{\prime}$ and also lies in a clique of cardinality $4.$ }

{\bf Proof.} The proof that $G^{\prime}$ is doubly homogeneously traceable is similar to that in the above proof of Lemma 1 (but easier).

Next suppose $v$ lies in a longest cycle of $G$ and $v$ lies in a clique of cardinality $4.$  Let $N(v)=\{x_1,x_2,x_3,x_4\}$ where $v,x_1,x_2,x_3$ form a clique and suppose $v$ is blown up into $K_4$ whose vertices are $v_1,v_2,v_3, v_4$ such that $v_i$ is adjacent to $x_i$ for $i=1,2,3,4.$ See Figure 2 for the change of local structures.
\vskip 3mm
\par
 \centerline{\includegraphics[width=3.3in]{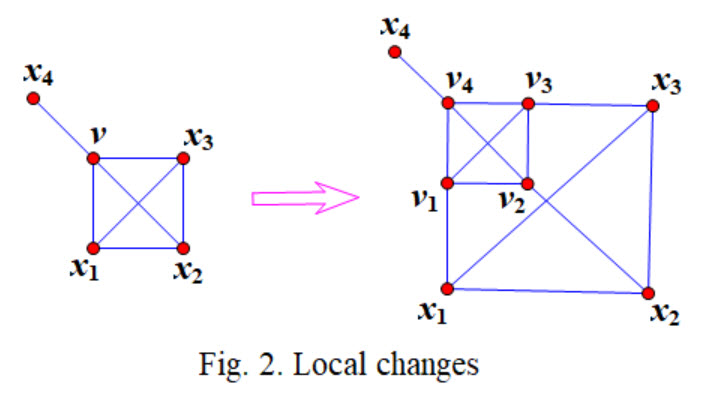}}
\par
Let $\ldots, x_i,v,x_j,\ldots$ be a longest cycle of $G$ with length $c.$ Then $G^{\prime}$ contains the cycle $\ldots, x_i,v_i,v_s,v_t,v_j,x_j,\ldots$
which has length $c+3.$ Thus ${\rm circum}(G^{\prime})\ge c+3.$ We then prove the reverse inequality. Let $C$ be a longest cycle of $G^{\prime},$  which has length at least $c+3.$ Denote $S=\{v_1,v_2,v_3,v_4\}.$ If $C$ contains no vertex from the set $S,$ it is also a cycle in $G$ and hence has length at most $c,$ a contradiction. Note that every vertex in $S$ has exactly one neighbor outside $S.$ Thus, if a cycle contains a vertex in $S,$ it contains at least two.
We have $|V(C)\cap S|\ge 2.$ If $w\in V(C)\cap S,$ then at least one neighbor of $w$ on $C$ belongs to $S.$
 Since $S$ is a clique and $C$ is a longest cycle, we deduce that  $|V(C)\cap S|=4.$ 
 On the cycle $C,$ using the vertex $v$ instead of $v_4$ and replacing a path of length $5$ by a path of length $2$, or replacing two paths of length $3$ and $2$ respectively by two edges we obtain a cycle of $G,$ where we have used the fact that $v,x_1,x_2,x_3$
 form a clique in $G.$ Since ${\rm circum}(G)=c,$ the cycle $C$ has length at most $c+3.$ This proves ${\rm circum}(G^{\prime})=c+3.$

 Finally we may choose $v_4$ as the vertex $v^{\prime}.$  \hfill $\Box$

 Now we are ready to state and prove the main results.

 {\bf Theorem 3.} {\it For every even integer $n\ge 10,$ there exists a cubic homogeneously traceable nonhamiltonian graph of order $n;$ for every integer $p\ge 18,$ there exists a $4$-regular homogeneously traceable graph of order $p$ and circumference $p-4.$ }

 {\bf Proof.} The Petersen graph $P$ depicted in Figure 3 is a cubic doubly homogeneously traceable graph of order $10$ and circumference $9.$
 \vskip 3mm
\par
 \centerline{\includegraphics[width=1.9in]{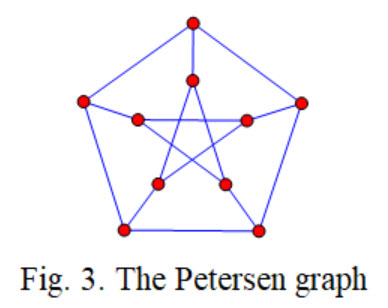}}
\par
Note that every vertex of $P$ lies in a longest cycle. Thus, choosing any vertex $v$ of $P$ and blowing up $v$ into $K_3$ we obtain a cubic graph $P_{12}$ of order $12.$ By Lemma 1, $P_{12}$ is doubly homogeneously traceable and has circumference $11.$ Let $u$ be a vertex of $P_{12}$ that lies in a longest cycle. In $P_{12},$
blowing up $u$ into $K_3$ we obtain a cubic graph $P_{14}$ of order $14.$ By Lemma 1, $P_{14}$ is doubly homogeneously traceable and has circumference $13.$ Continuing this process we can construct a cubic homogeneously traceable graph of order $n$ and circumference $n-1$ for any even integer $n\ge 10.$

It is easy to verify that the three graphs in Figures 4, 5 and 6 are $4$-regular doubly homogeneously traceable graphs of order $p$ and circumference $p-4$ for $p=18,19,20$
respectively, where the vertices $x,y,z$ lie in a longest cycle and in a clique of cardinality $4.$
\vskip 3mm
\par
 \centerline{\includegraphics[width=3.3in]{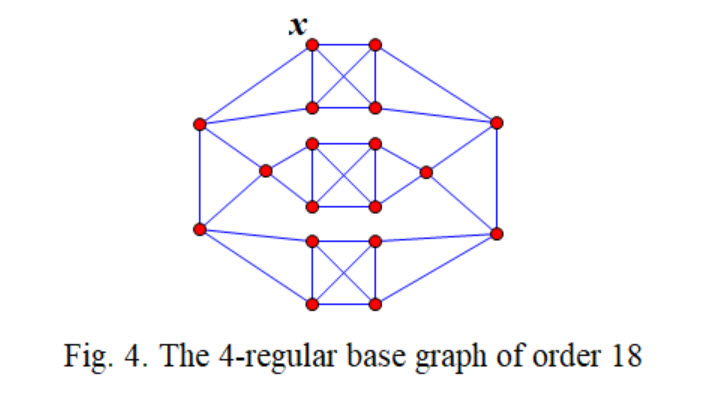}}
\par
\vskip 3mm
\par
 \centerline{\includegraphics[width=3in]{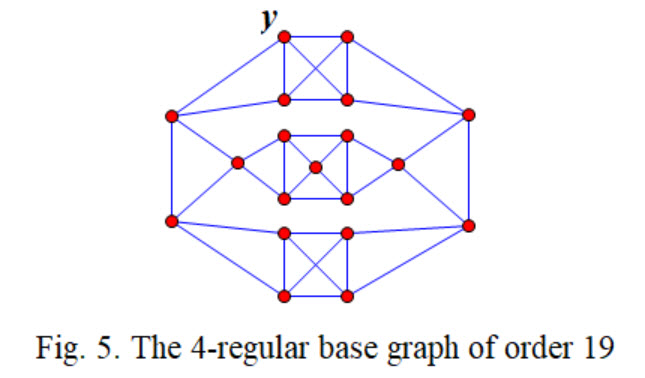}}
\par
\vskip 3mm
\par
 \centerline{\includegraphics[width=2.9in]{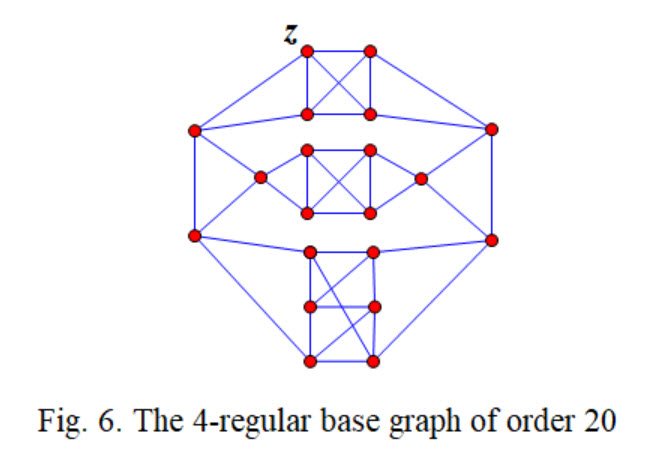}}
\par
Next we apply Lemma 2 repeatedly. Starting from the graph in Figure 4 and the vertex $x$, successively blowing up a vertex that lies in a longest cycle and in a clique of cardinality $4,$ we can construct a $4$-regular homogeneously traceable graph of order $p$ and circumference $p-4$ for every integer $p\ge 18$ with $p\equiv 0\,\,\,{\rm mod}\,\,3.$ Starting from the graph in Figure 5 and the vertex $y$, successively blowing up a vertex that lies in a longest cycle and in a clique of cardinality $4,$ we can construct a $4$-regular homogeneously traceable graph of order $p$ and circumference $p-4$ for every integer $p\ge 19$ with $p\equiv 1\,\,\,{\rm mod}\,\,3.$ Starting from the graph in Figure 6 and the vertex $z$, successively blowing up a vertex that lies in a longest cycle and in a clique of cardinality $4,$ we can construct a $4$-regular homogeneously traceable graph of order $p$ and circumference $p-4$ for every integer $p\ge 20$ with $p\equiv 2\,\,\,{\rm mod}\,\,3.$
This completes the proof. \hfill $\Box$

{\bf Remark.} The above proof of Theorem 3 shows that in the statement of Theorem 3, we may replace ``homogeneously traceable" by ``doubly homogeneously traceable". But we prefer
the current version, since the term ``doubly homogeneously traceable" is technical in some sense.

\section{Unsolved problems}

It is known ([2, Theorem 4] and [4, pp. 9-11]) that the minimum size of a homogeneously traceable nonhamiltonian graph of order $n$ is
$\lceil 5n/4\rceil.$

The extremal problem concerning the independence number is easy.

{\bf Theorem 4.} {\it The maximum independence number of a homogeneously traceable graph of order $n$ is $\lfloor n/2\rfloor.$}

{\bf Proof.} Let $G$ be a homogeneously traceable graph of order $n,$ and let $v_1,v_2,\ldots,v_n$ be a Hamilton path.
Suppose $S$ is an independent set of $G.$ If $n$ is even, $S$ contains at most one vertex in each of the edges $v_1v_2,$ $v_3v_4,$ $\ldots,$ $v_{n-1}v_n$
and hence $|S|\le n/2.$ Now suppose $n$ is odd. Similarly, we have $|S|\le (n+1)/2.$ We will show that $|S|$ cannot equal $(n+1)/2.$ To the contrary, assume
$|S|=(n+1)/2.$ Then $S=\{v_1,v_3,v_5,\ldots, v_n\}.$ Since $G$ is homogeneously traceable, there is a Hamilton path $v_2,v_{i_2},\ldots,v_{i_n}.$
Since $n$ is odd, there exists an integer $k$ with $2\le k\le n-1$ such that both $i_k$ and $i_{k+1}$ are odd. But $v_{i_k}$ and $v_{i_{k+1}}$ are adjacent
and they both belong to $S,$ contradicting the condition that $S$ is an independent set. It follows that $|S|\le (n-1)/2.$ We have proved that
$|S|\le \lfloor n/2\rfloor.$

This upper bound $\lfloor n/2\rfloor$ can be attained by the cycle $C_n$ which is homogeneously traceable, and hence it is indeed the maximum value. \hfill $\Box$

Finally we pose two unsolved problems.

{\bf Conjecture 1.} {\it The minimum circumference of a homogeneously traceable graph of order $n$ is $\lceil 2n/3\rceil+2.$}

The circumference $\lceil 2n/3\rceil+2$ in Conjecture 1 is attained by the graph in Figure 7 where $p=\lfloor (n-6)/3\rfloor$ and 
when $p\ge 2$ the vertices $u$ and $v$ are distinct, $x$ and $y$ are distinct and $w$ and $z$ are distinct.
\vskip 3mm
\par
 \centerline{\includegraphics[width=3.6in]{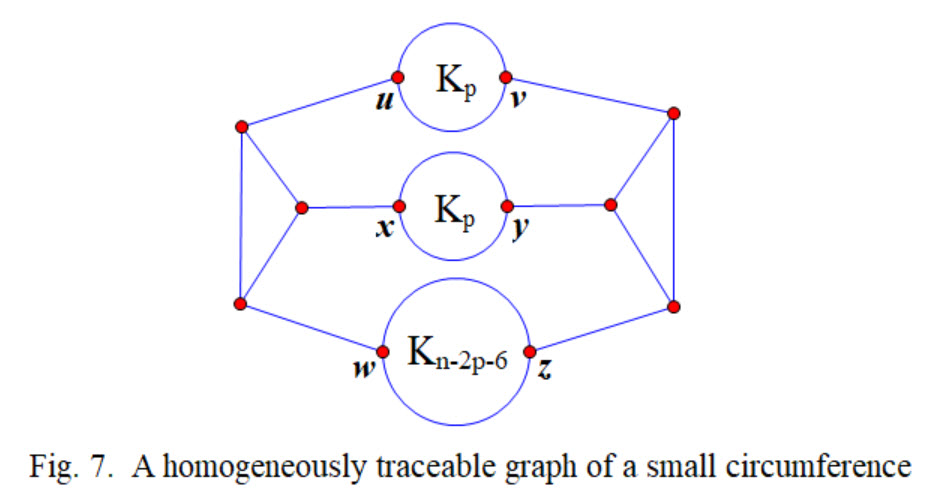}}
\par

{\bf Problem 2.} {\it Given an integer $k\ge 4,$ determine the integers $n$ such that there exists a $k$-regular homogeneously traceable nonhamiltonian 
graph of order $n.$}

Theorem 3 solves the case $k=3$ of Problem 2. A computer search shows that there exists no $4$-regular homogeneously traceable nonhamiltonian
graph of order $\le 15.$ Thus, according to Theorem 3, only the two orders $16$ and $17$ are uncertain for $k=4.$

\vskip 5mm
{\bf Acknowledgement.} This research  was supported by the NSFC grants 11671148 and 11771148 and Science and Technology Commission of Shanghai Municipality (STCSM) grant 18dz2271000.

\end{document}